# MEASURING AND TESTING DEPENDENCE BY CORRELATION OF DISTANCES

By Gábor J. Székely,[1] Maria L. Rizzo and Nail K. Bakirov

*Bowling Green State University, Bowling Green State University and USC Russian Academy of Sciences*

Distance correlation is a new measure of dependence between random vectors. Distance covariance and distance correlation are analogous to product-moment covariance and correlation, but unlike the classical definition of correlation, distance correlation is zero only if the random vectors are independent. The empirical distance dependence measures are based on certain Euclidean distances between sample elements rather than sample moments, yet have a compact representation analogous to the classical covariance and correlation. Asymptotic properties and applications in testing independence are discussed. Implementation of the test and Monte Carlo results are also presented.

**1. Introduction.** Distance correlation provides a new approach to the problem of testing the joint independence of random vectors. For all distributions with finite first moments, distance correlation $\mathcal{R}$ generalizes the idea of correlation in two fundamental ways:

(i) $\mathcal{R}(X, Y)$ is defined for $X$ and $Y$ in arbitrary dimensions;
(ii) $\mathcal{R}(X, Y) = 0$ characterizes independence of $X$ and $Y$.

Distance correlation has properties of a true dependence measure, analogous to product-moment correlation $\rho$. Distance correlation satisfies $0 \leq \mathcal{R} \leq 1$, and $\mathcal{R} = 0$ only if $X$ and $Y$ are independent. In the bivariate normal case, $\mathcal{R}$ is a function of $\rho$, and $\mathcal{R}(X, Y) \leq |\rho(X, Y)|$ with equality when $\rho = \pm 1$.

Throughout this paper $X$ in $\mathbb{R}^p$ and $Y$ in $\mathbb{R}^q$ are random vectors, where $p$ and $q$ are positive integers. The characteristic functions of $X$ and $Y$ are denoted $f_X$ and $f_Y$, respectively, and the joint characteristic function of $X$

Received November 2006; revised March 2007.
[1]Supported by the National Science Foundation while working at the Foundation.

and $Y$ is denoted $f_{X,Y}$. Distance covariance $\mathcal{V}$ can be applied to measure the distance $\|f_{X,Y}(t,s) - f_X(t)f_Y(s)\|$ between the joint characteristic function and the product of the marginal characteristic functions (the norm $\|\cdot\|$ is defined in Section 2), and to test the hypothesis of independence

$$H_0 : f_{X,Y} = f_X f_Y \quad \text{vs.} \quad H_1 : f_{X,Y} \neq f_X f_Y.$$

The importance of the independence assumption for inference arises, for example, in clinical studies with the case-only design, which uses only diseased subjects assumed to be independent in the study population. In this design, inferences on multiplicative gene interactions (see [1]) can be highly distorted when there is a departure from independence. Classical methods such as the Wilks Lambda [14] or Puri–Sen [8] likelihood ratio tests are not applicable if the dimension exceeds the sample size, or when distributional assumptions do not hold (see, e.g., [7] regarding the prevalence of nonnormality in biology and ecology). A further limitation of multivariate extensions of methods based on ranks is that they are ineffective for testing nonmonotone types of dependence.

We propose an omnibus test of independence that is easily implemented in arbitrary dimension. In our Monte Carlo results the distance covariance test exhibits superior power against nonmonotone types of dependence while maintaining good power performance in the multivariate normal case relative to the parametric likelihood ratio test. Distance correlation can also be applied as an index of dependence; for example, in meta-analysis [12] distance correlation would be a more generally applicable index than product-moment correlation, without requiring normality for valid inferences.

Theoretical properties of distance covariance and correlation are covered in Section 2, extensions in Section 3, and results for the bivariate normal case in Section 4. Empirical results are presented in Section 5, followed by a summary in Section 6.

## 2. Theoretical properties of distance dependence measures.

*Notation.* The scalar product of vectors $t$ and $s$ is denoted by $\langle t, s \rangle$. For complex-valued functions $f(\cdot)$, the complex conjugate of $f$ is denoted by $\overline{f}$ and $|f|^2 = f\overline{f}$. The Euclidean norm of $x$ in $\mathbb{R}^p$ is $|x|_p$. A sample from the distribution of $X$ in $\mathbb{R}^p$ is denoted by the $n \times p$ matrix $\mathbf{X}$, and the sample vectors (rows) are labeled $X_1, \ldots, X_n$. A primed variable $X'$ is an independent copy of $X$; that is, $X$ and $X'$ are independent and identically distributed.

DEFINITION 1. For complex functions $\gamma$ defined on $\mathbb{R}^p \times \mathbb{R}^q$ the $\|\cdot\|_w$-norm in the weighted $L_2$ space of functions on $\mathbb{R}^{p+q}$ is defined by

$$(2.1) \qquad \|\gamma(t,s)\|_w^2 = \int_{\mathbb{R}^{p+q}} |\gamma(t,s)|^2 w(t,s)\, dt\, ds,$$

where $w(t,s)$ is an arbitrary positive weight function for which the integral above exists.

2.1. *Choice of weight function.* Using the $\|\cdot\|_w$-norm (2.1) with a suitable choice of weight function $w(t,s)$, we define a measure of dependence

$$\begin{aligned}(2.2)\qquad \mathcal{V}^2(X,Y;w) &= \|f_{X,Y}(t,s) - f_X(t)f_Y(s)\|_w^2 \\ &= \int_{\mathbb{R}^{p+q}} |f_{X,Y}(t,s) - f_X(t)f_Y(s)|^2 w(t,s)\,dt\,ds,\end{aligned}$$

such that $\mathcal{V}^2(X,Y;w) = 0$ if and only if $X$ and $Y$ are independent. In this paper $\mathcal{V}$ will be analogous to the absolute value of the classical product-moment covariance. If we divide $\mathcal{V}(X,Y;w)$ by $\sqrt{\mathcal{V}(X;w)\mathcal{V}(Y;w)}$, where

$$\mathcal{V}^2(X;w) = \int_{\mathbb{R}^{2p}} |f_{X,X}(t,s) - f_X(t)f_X(s)|^2 w(t,s)\,dt\,ds,$$

we have a type of unsigned correlation $\mathcal{R}_w$.

Not every weight function leads to an "interesting" $\mathcal{R}_w$, however. The coefficient $\mathcal{R}_w$ should be scale invariant, that is, invariant with respect to transformations $(X,Y) \mapsto (\epsilon X, \epsilon Y)$, for $\epsilon > 0$. We also require that $\mathcal{R}_w$ is positive for dependent variables. It is easy to check that if the weight function $w(t,s)$ is *integrable*, and both $X$ and $Y$ have finite variance, then the Taylor expansions of the underlying characteristic functions show that

$$\lim_{\epsilon \to 0} \frac{\mathcal{V}^2(\epsilon X, \epsilon Y; w)}{\sqrt{\mathcal{V}^2(\epsilon X; w)\mathcal{V}^2(\epsilon Y; w)}} = \rho^2(X,Y).$$

Thus for integrable $w$, if $\rho = 0$, then $\mathcal{R}_w$ can be arbitrarily close to zero even if $X$ and $Y$ are dependent. However, by applying a nonintegrable weight function, we obtain an $\mathcal{R}_w$ that is scale invariant and cannot be zero for dependent $X$ and $Y$. We do not claim that our choice for $w$ is the only reasonable one, but it will become clear in the following sections that our choice (2.4) results in very simple and applicable empirical formulas. (A more complicated weight function is applied in [2], which leads to a more computationally difficult statistic and does not have the interesting correlation form.)

The crucial observation is the following lemma.

LEMMA 1. *If $0 < \alpha < 2$, then for all $x$ in $\mathbb{R}^d$*

$$\int_{\mathbb{R}^d} \frac{1 - \cos\langle t, x \rangle}{|t|_d^{d+\alpha}}\,dt = C(d,\alpha)|x|^\alpha,$$

*where*

$$C(d,\alpha) = \frac{2\pi^{d/2}\Gamma(1-\alpha/2)}{\alpha 2^\alpha \Gamma((d+\alpha)/2)}$$

and $\Gamma(\cdot)$ is the complete gamma function. The integrals at $0$ and $\infty$ are meant in the principal value sense: $\lim_{\varepsilon \to 0} \int_{\mathbb{R}^d \setminus \{\varepsilon B + \varepsilon^{-1} B^c\}}$, where $B$ is the unit ball (centered at 0) in $\mathbb{R}^d$ and $B^c$ is the complement of $B$.

See [11] for the proof of Lemma 1. In the simplest case, $\alpha = 1$, the constant in Lemma 1 is

$$c_d = C(d, 1) = \frac{\pi^{(1+d)/2}}{\Gamma((1+d)/2)}. \tag{2.3}$$

In view of Lemma 1, it is natural to choose the weight function

$$w(t, s) = (c_p c_q |t|_p^{1+p} |s|_q^{1+q})^{-1}, \tag{2.4}$$

corresponding to $\alpha = 1$. We apply the weight function (2.4) and the corresponding weighted $L_2$ norm $\|\cdot\|$, omitting the index $w$, and write the dependence measure (2.2) as $\mathcal{V}^2(X, Y)$. In integrals we also use the symbol $d\omega$, which is defined by

$$d\omega = (c_p c_q |t|_p^{1+p} |s|_q^{1+q})^{-1} \, dt \, ds.$$

For finiteness of $\|f_{X,Y}(t, s) - f_X(t) f_Y(s)\|^2$ it is sufficient that $E|X|_p < \infty$ and $E|Y|_q < \infty$. By the Cauchy–Bunyakovsky inequality

$$|f_{X,Y}(t,s) - f_X(t)f_Y(s)|^2 = [E(e^{i\langle t,X\rangle} - f_X(t))(e^{i\langle s,Y\rangle} - f_Y(s))]^2$$
$$\leq E[e^{i\langle t,X\rangle} - f_X(t)]^2 E[e^{i\langle s,Y\rangle} - f_Y(s)]^2$$
$$= (1 - |f_X(t)|^2)(1 - |f_Y(s)|^2).$$

If $E(|X|_p + |Y|_q) < \infty$, then by Lemma 1 and by Fubini's theorem it follows that

$$\int_{\mathbb{R}^{p+q}} |f_{X,Y}(t, s) - f_X(t) f_Y(s)|^2 \, d\omega$$
$$\leq \int_{\mathbb{R}^p} \frac{1 - |f_X(t)|^2}{c_p |t|_p^{1+p}} \, dt \int_{\mathbb{R}^q} \frac{1 - |f_Y(s)|^2}{c_q |s|_q^{1+q}} \, ds$$
$$= E\left[\int_{\mathbb{R}^p} \frac{1 - \cos\langle t, X - X'\rangle}{c_p |t|_p^{1+p}} \, dt\right] \cdot E\left[\int_{\mathbb{R}^q} \frac{1 - \cos\langle s, Y - Y'\rangle}{c_q |s|_q^{1+q}} \, ds\right]$$
$$= E|X - X'|_p E|Y - Y'|_q < \infty. \tag{2.5}$$

Thus we have the following definitions.

### 2.2. *Distance covariance and distance correlation.*

DEFINITION 2 (*Distance covariance*). The distance covariance (dCov) between random vectors $X$ and $Y$ with finite first moments is the nonnegative number $\mathcal{V}(X,Y)$ defined by

$$
\begin{aligned}
\mathcal{V}^2(X,Y) &= \|f_{X,Y}(t,s) - f_X(t)f_Y(s)\|^2 \\
&= \frac{1}{c_p c_q} \int_{\mathbb{R}^{p+q}} \frac{|f_{X,Y}(t,s) - f_X(t)f_Y(s)|^2}{|t|_p^{1+p}|s|_q^{1+q}} \, dt \, ds.
\end{aligned}
$$
(2.6)

Similarly, distance variance (dVar) is defined as the square root of

$$\mathcal{V}^2(X) = \mathcal{V}^2(X,X) = \|f_{X,X}(t,s) - f_X(t)f_X(s)\|^2.$$

REMARK 1. If $E(|X|_p + |Y|_q) = \infty$ but $E(|X|_p^\alpha + |Y|_q^\alpha) < \infty$ for some $0 < \alpha < 1$, then one can apply $\mathcal{V}^{(\alpha)}$ and $\mathcal{R}^{(\alpha)}$ (see Section 3.1); otherwise one can apply a suitable transformation of $(X,Y)$ into bounded random variables $(\widetilde{X}, \widetilde{Y})$ such that $\widetilde{X}$ and $\widetilde{Y}$ are independent if and only if $X$ and $Y$ are independent.

DEFINITION 3 (*Distance correlation*). The distance correlation (dCor) between random vectors $X$ and $Y$ with finite first moments is the nonnegative number $\mathcal{R}(X,Y)$ defined by

$$
(2.7) \quad \mathcal{R}^2(X,Y) = \begin{cases} \dfrac{\mathcal{V}^2(X,Y)}{\sqrt{\mathcal{V}^2(X)\mathcal{V}^2(Y)}}, & \mathcal{V}^2(X)\mathcal{V}^2(Y) > 0, \\ 0, & \mathcal{V}^2(X)\mathcal{V}^2(Y) = 0. \end{cases}
$$

Clearly the definition of $\mathcal{R}$ in (2.7) suggests an analogy with the product-moment correlation coefficient $\rho$. Analogous properties are established in Theorem 3. The relation between $\mathcal{V}$, $\mathcal{R}$ and $\rho$ in the bivariate normal case will be established in Theorem 7.

The distance dependence statistics are defined as follows. For an observed random sample $(\mathbf{X}, \mathbf{Y}) = \{(X_k, Y_k) : k = 1, \dots, n\}$ from the joint distribution of random vectors $X$ in $\mathbb{R}^p$ and $Y$ in $\mathbb{R}^q$, define

$$a_{kl} = |X_k - X_l|_p, \qquad \bar{a}_{k\cdot} = \frac{1}{n}\sum_{l=1}^n a_{kl}, \qquad \bar{a}_{\cdot l} = \frac{1}{n}\sum_{k=1}^n a_{kl},$$

$$\bar{a}_{\cdot\cdot} = \frac{1}{n^2}\sum_{k,l=1}^n a_{kl}, \qquad A_{kl} = a_{kl} - \bar{a}_{k\cdot} - \bar{a}_{\cdot l} + \bar{a}_{\cdot\cdot},$$

$k, l = 1, \ldots, n$. Similarly, define $b_{kl} = |Y_k - Y_l|_q$ and $B_{kl} = b_{kl} - \bar{b}_{k\cdot} - \bar{b}_{\cdot l} + \bar{b}_{\cdot\cdot}$, for $k, l = 1, \ldots, n$.

DEFINITION 4. The empirical distance covariance $\mathcal{V}_n(\mathbf{X}, \mathbf{Y})$ is the nonnegative number defined by

$$\mathcal{V}_n^2(\mathbf{X}, \mathbf{Y}) = \frac{1}{n^2} \sum_{k,l=1}^n A_{kl} B_{kl}. \tag{2.8}$$

Similarly, $\mathcal{V}_n(\mathbf{X})$ is the nonnegative number defined by

$$\mathcal{V}_n^2(\mathbf{X}) = \mathcal{V}_n^2(\mathbf{X}, \mathbf{X}) = \frac{1}{n^2} \sum_{k,l=1}^n A_{kl}^2. \tag{2.9}$$

Although it may not be immediately obvious that $\mathcal{V}_n^2(\mathbf{X}, \mathbf{Y}) \geq 0$, this fact as well as the motivation for the definition of $\mathcal{V}_n$ will be clear from Theorem 1 below.

DEFINITION 5. The empirical distance correlation $\mathcal{R}_n(\mathbf{X}, \mathbf{Y})$ is the square root of

$$\mathcal{R}_n^2(\mathbf{X}, \mathbf{Y}) = \begin{cases} \dfrac{\mathcal{V}_n^2(\mathbf{X}, \mathbf{Y})}{\sqrt{\mathcal{V}_n^2(\mathbf{X})\mathcal{V}_n^2(\mathbf{Y})}}, & \mathcal{V}_n^2(\mathbf{X})\mathcal{V}_n^2(\mathbf{Y}) > 0, \\ 0, & \mathcal{V}_n^2(\mathbf{X})\mathcal{V}_n^2(\mathbf{Y}) = 0. \end{cases} \tag{2.10}$$

REMARK 2. The statistic $\mathcal{V}_n(\mathbf{X}) = 0$ if and only if every sample observation is identical. Indeed, if $\mathcal{V}_n(\mathbf{X}) = 0$, then $A_{kl} = 0$ for $k, l = 1, \ldots, n$. Thus $0 = A_{kk} = -\bar{a}_{k\cdot} - \bar{a}_{\cdot k} + \bar{a}_{\cdot\cdot}$ implies that $\bar{a}_{k\cdot} = \bar{a}_{\cdot k} = \bar{a}_{\cdot\cdot}/2$, and

$$0 = A_{kl} = a_{kl} - \bar{a}_{k\cdot} - \bar{a}_{\cdot l} + \bar{a}_{\cdot\cdot} = a_{kl} = |X_k - X_l|_p,$$

so $X_1 = \cdots = X_n$.

It is clear that $\mathcal{R}_n$ is easy to compute, and in the following sections it will be shown that $\mathcal{R}_n$ is a good empirical measure of dependence.

2.3. *Properties of distance covariance.* It would have been natural, but less elementary, to define $\mathcal{V}_n(\mathbf{X}, \mathbf{Y})$ as $\|f_{X,Y}^n(t, s) - f_X^n(t) f_Y^n(s)\|$, where

$$f_{X,Y}^n(t, s) = \frac{1}{n} \sum_{k=1}^n \exp\{i \langle t, X_k \rangle + i \langle s, Y_k \rangle\}$$

is the empirical characteristic function of the sample, $\{(X_1, Y_1), \ldots, (X_n, Y_n)\}$, and

$$f_X^n(t) = \frac{1}{n} \sum_{k=1}^n \exp\{i \langle t, X_k \rangle\}, \qquad f_Y^n(s) = \frac{1}{n} \sum_{k=1}^n \exp\{i \langle s, Y_k \rangle\}$$

are the marginal empirical characteristic functions of the $X$ sample and $Y$ sample, respectively. Our first theorem shows that the two definitions are equivalent.

THEOREM 1. *If $(\mathbf{X}, \mathbf{Y})$ is a sample from the joint distribution of $(X, Y)$, then*

$$\mathcal{V}_n^2(\mathbf{X}, \mathbf{Y}) = \|f_{X,Y}^n(t,s) - f_X^n(t)f_Y^n(s)\|^2.$$

PROOF. Lemma 1 implies that there exist constants $c_p$ and $c_q$ such that for all $X$ in $\mathbb{R}^p$, $Y$ in $\mathbb{R}^q$,

$$(2.11) \qquad \int_{\mathbb{R}^p} \frac{1 - \exp\{i\langle t, X\rangle\}}{|t|_p^{1+p}}\, dt = c_p |X|_p,$$

$$(2.12) \qquad \int_{\mathbb{R}^q} \frac{1 - \exp\{i\langle s, Y\rangle\}}{|s|_q^{1+q}}\, ds = c_q |Y|_q,$$

$$(2.13) \qquad \int_{\mathbb{R}^p}\int_{\mathbb{R}^q} \frac{1 - \exp\{i\langle t, X\rangle + i\langle s, Y\rangle\}}{|t|_p^{1+p}|s|_q^{1+q}}\, dt\, ds = c_p c_q |X|_p |Y|_q,$$

where the integrals are understood in the principal value sense. For simplicity, consider the case $p = q = 1$. The distance between the empirical characteristic functions in the weighted norm $w(t,s) = \pi^{-2}t^{-2}s^{-2}$ involves $|f_{X,Y}^n(t,s)|^2$, $|f_X^n(t)f_Y^n(s)|^2$ and $\overline{f_{X,Y}^n(t,s)}f_X^n(t)f_Y^n(s)$. For the first we have

$$f_{X,Y}^n(t,s)\overline{f_{X,Y}^n(t,s)} = \frac{1}{n^2}\sum_{k,l=1}^n \cos(X_k - X_l)t\cos(Y_k - Y_l)s + V_1,$$

where $V_1$ represents terms that vanish when the integral $\|f_{X,Y}^n(t,s) - f_X^n(t)f_Y^n(s)\|^2$ is evaluated. The second expression is

$$f_X^n(t)f_Y^n(s)\overline{f_X^n(t)f_Y^n(s)}$$
$$= \frac{1}{n^2}\sum_{k,l=1}^n \cos(X_k - X_l)t \frac{1}{n^2}\sum_{k,l=1}^n \cos(Y_k - Y_l)s + V_2$$

and the third is

$$f_{X,Y}^n(t,s)\overline{f_X^n(t)f_Y^n(s)}$$
$$= \frac{1}{n^3}\sum_{k,l,m=1}^n \cos(X_k - X_l)t\cos(Y_k - Y_m)s + V_3,$$

where $V_2$ and $V_3$ represent terms that vanish when the integral is evaluated. To evaluate the integral $\|f_{X,Y}^n(t,s) - f_X^n(t)f_Y^n(s)\|^2$, apply Lemma 1, and statements (2.11), (2.12) and (2.13) using

$$\cos u \cos v = 1 - (1 - \cos u) - (1 - \cos v) + (1 - \cos u)(1 - \cos v).$$

After cancellation in the numerator of the integrand it remains to evaluate integrals of the type

$$\int_{\mathbb{R}^2} (1 - \cos(X_k - X_l)t)(1 - \cos(Y_k - Y_l)s) \frac{dt}{t^2} \frac{ds}{s^2}$$

$$= \int_{\mathbb{R}} (1 - \cos(X_k - X_l)t) \frac{dt}{t^2}$$

$$\times \int_{\mathbb{R}} (1 - \cos(Y_k - Y_l)s) \frac{ds}{s^2}$$

$$= c_1^2 |X_k - X_l| \, |Y_k - Y_l|.$$

For random vectors $X$ in $\mathbb{R}^p$ and $Y$ in $\mathbb{R}^q$, the same steps are applied using $w(t,s) = \{c_p c_q |t|_p^{1+p} |s|_q^{1+q}\}^{-1}$. Thus

$$\|f_{X,Y}^n(t,s) - f_X^n(t) f_Y^n(s)\|^2 = S_1 + S_2 - 2S_3, \tag{2.14}$$

where

$$S_1 = \frac{1}{n^2} \sum_{k,l=1}^n |X_k - X_l|_p |Y_k - Y_l|_q, \tag{2.15}$$

$$S_2 = \frac{1}{n^2} \sum_{k,l=1}^n |X_k - X_l|_p \frac{1}{n^2} \sum_{k,l=1}^n |Y_k - Y_l|_q, \tag{2.16}$$

$$S_3 = \frac{1}{n^3} \sum_{k=1}^n \sum_{l,m=1}^n |X_k - X_l|_p |Y_k - Y_m|_q. \tag{2.17}$$

To complete the proof we need to verify the algebraic identity

$$\mathcal{V}_n^2(\mathbf{X}, \mathbf{Y}) = S_1 + S_2 - 2S_3. \tag{2.18}$$

For the proof of (2.18) see the Appendix. Then (2.14) and (2.18) imply that $\mathcal{V}_n^2(\mathbf{X}, \mathbf{Y}) = \|f_{X,Y}^n(t,s) - f_X^n(t) f_Y^n(s)\|^2$. □

THEOREM 2. *If* $E|X|_p < \infty$ *and* $E|Y|_q < \infty$, *then almost surely*

$$\lim_{n \to \infty} \mathcal{V}_n(\mathbf{X}, \mathbf{Y}) = \mathcal{V}(X, Y). \tag{2.19}$$

PROOF. Define

$$\xi_n(t,s) = \frac{1}{n} \sum_{k=1}^n e^{i\langle t, X_k \rangle + i\langle s, Y_k \rangle} - \frac{1}{n} \sum_{k=1}^n e^{i\langle t, X_k \rangle} \frac{1}{n} \sum_{k=1}^n e^{i\langle s, Y_k \rangle},$$

so that $\mathcal{V}_n^2 = \|\xi_n(t,s)\|^2$. Then after elementary transformations

$$\xi_n(t,s) = \frac{1}{n} \sum_{k=1}^n u_k v_k - \frac{1}{n} \sum_{k=1}^n u_k \frac{1}{n} \sum_{k=1}^n v_k,$$

where $u_k = \exp(i\langle t, X_k\rangle) - f_X(t)$ and $v_k = \exp(i\langle s, Y_k\rangle) - f_Y(s)$.

For each $\delta > 0$ define the region

(2.20) $$D(\delta) = \{(t,s) : \delta \leq |t|_p \leq 1/\delta, \delta \leq |s|_q \leq 1/\delta\}$$

and random variables

$$\mathcal{V}_{n,\delta}^2 = \int_{D(\delta)} |\xi_n(t,s)|^2 \, d\omega.$$

For any fixed $\delta > 0$, the weight function $w(t,s)$ is bounded on $D(\delta)$. Hence $\mathcal{V}_{n,\delta}^2$ is a combination of $V$-statistics of bounded random variables. For each $\delta > 0$ by the strong law of large numbers (SLLN) for $V$-statistics, it follows that almost surely

$$\lim_{n\to\infty} \mathcal{V}_{n,\delta}^2 = \mathcal{V}_{\cdot,\delta}^2 = \int_{D(\delta)} |f_{X,Y}(t,s) - f_X(t)f_Y(s)|^2 \, d\omega.$$

Clearly $\mathcal{V}_{\cdot,\delta}^2$ converges to $\mathcal{V}^2$ as $\delta$ tends to zero. Now it remains to prove that almost surely

(2.21) $$\limsup_{\delta\to 0} \limsup_{n\to\infty} |\mathcal{V}_{n,\delta}^2 - \mathcal{V}_n^2| = 0.$$

For each $\delta > 0$

(2.22) $$|\mathcal{V}_{n,\delta}^2 - \mathcal{V}_n^2| \leq \int_{|t|_p < \delta} |\xi_n(t,s)|^2 \, d\omega + \int_{|t|_p > 1/\delta} |\xi_n(t,s)|^2 \, d\omega$$
$$+ \int_{|s|_q < \delta} |\xi_n(t,s)|^2 \, d\omega + \int_{|s|_q > 1/\delta} |\xi_n(t,s)|^2 \, d\omega.$$

For $z = (z_1, z_2, \ldots, z_p)$ in $\mathbb{R}^p$ define the function

$$G(y) = \int_{|z| < y} \frac{1 - \cos z_1}{|z|^{1+p}} \, dz.$$

Clearly $G(y)$ is bounded by $c_p$ and $\lim_{y\to 0} G(y) = 0$. Applying the inequality $|x+y|^2 \leq 2|x|^2 + 2|y|^2$ and the Cauchy–Bunyakovsky inequality for sums, one can obtain that

(2.23) $$|\xi_n(t,s)|^2 \leq 2\left|\frac{1}{n}\sum_{k=1}^n u_k v_k\right|^2 + 2\left|\frac{1}{n}\sum_{k=1}^n u_k \frac{1}{n}\sum_{k=1}^n v_k\right|^2$$
$$\leq \frac{4}{n}\sum_{k=1}^n |u_k|^2 \frac{1}{n}\sum_{k=1}^n |v_k|^2.$$

Hence the first summand in (2.22) satisfies

(2.24) $$\int_{|t|_p < \delta} |\xi_n(t,s)|^2 \, d\omega \leq \frac{4}{n}\sum_{k=1}^n \int_{|t|_p < \delta} \frac{|u_k|^2 \, dt}{c_p |t|_p^{1+p}} \frac{1}{n}\sum_{k=1}^n \int_{\mathbb{R}^q} \frac{|v_k|^2 \, ds}{c_q |s|_q^{1+q}}.$$

Here $|v_k|^2 = 1 + |f_Y(s)|^2 - e^{i\langle s, Y_k\rangle}\overline{f_Y(s)} - e^{-i\langle s, Y_k\rangle}f_Y(s)$, thus

$$\int_{\mathbb{R}^q} \frac{|v_k|^2\, ds}{c_q |s|_q^{1+q}} = (2E_Y|Y_k - Y| - E|Y - Y'|) \leq 2(|Y_k| + E|Y|),$$

where the expectation $E_Y$ is taken with respect to $Y$, and $Y' \stackrel{D}{=} Y$ is independent of $Y_k$. Further, after a suitable change of variables

$$\int_{|t|_p < \delta} \frac{|u_k|^2\, dt}{c_p |t|_p^{1+p}} = 2E_X|X_k - X|G(|X_k - X|\delta) - E|X - X'|G(|X - X'|\delta)$$

$$\leq 2E_X|X_k - X|G(|X_k - X|\delta),$$

where the expectation $E_X$ is taken with respect to $X$, and $X' \stackrel{D}{=} X$ is independent of $X_k$. Therefore, from (2.24)

$$\int_{|t|_p < \delta} |\xi_n(t,s)|^2\, d\omega$$

$$\leq 4\frac{2}{n}\sum_{k=1}^n (|Y_k| + E|Y|)\frac{2}{n}\sum_{k=1}^n E_X|X_k - X|G(|X_k - X|\delta).$$

By the SLLN

$$\limsup_{n\to\infty} \int_{|t|_p < \delta} |\xi_n(t,s)|^2\, d\omega \leq 4 \cdot 2 \cdot 2E|Y| \cdot 2E|X_1 - X_2|G(|X_1 - X_2|\delta)$$

almost surely. Therefore by the Lebesgue bounded convergence theorem for integrals and expectations

$$\limsup_{\delta\to 0}\limsup_{n\to\infty} \int_{|t|_p < \delta} |\xi_n(t,s)|^2\, d\omega = 0$$

almost surely.

Consider now the second summand in (2.22). Inequalities (2.23) imply that $|u_k|^2 \leq 4$ and $\frac{1}{n}\sum_{k=1}^n |u_k|^2 \leq 4$, hence

$$\int_{|t|_p > 1/\delta} \frac{|u_k|^2\, dt}{c_p |t|_p^{1+p}} \leq 16 \int_{|t|_p > 1/\delta} \frac{dt}{c_p |t|_p^{1+p}} \int_{\mathbb{R}^q} \frac{1}{n}\sum_{k=1}^n |v_k|^2 \frac{ds}{c_q |s|_q^{1+q}}$$

$$\leq 16\delta \frac{2}{n}\sum_{k=1}^n (|Y_k| + E|Y|).$$

Thus, almost surely

$$\limsup_{\delta\to 0}\limsup_{n\to\infty} \int_{|t|_p > 1/\delta} |\xi_n(t,s)|^2\, d\omega = 0.$$

One can apply a similar argument to the remaining summands in (2.22) to obtain (2.21). $\square$

COROLLARY 1. *If $E(|X|_p + |Y|_q) < \infty$, then almost surely,*
$$\lim_{n\to\infty} \mathcal{R}_n^2(\mathbf{X},\mathbf{Y}) = \mathcal{R}^2(X,Y).$$

The definition of dCor suggests that our distance dependence measures are analogous in at least some respects to the corresponding product-moment correlation. By analogy, certain properties of classical correlation and variance definitions should also hold for dCor and dVar. These properties are established in Theorems 3 and 4.

THEOREM 3 (Properties of dCor).

(i) *If $E(|X|_p + |Y|_q) < \infty$, then $0 \leq \mathcal{R} \leq 1$, and $\mathcal{R}(X,Y) = 0$ if and only if $X$ and $Y$ are independent.*

(ii) $0 \leq \mathcal{R}_n \leq 1$.

(iii) *If $\mathcal{R}_n(\mathbf{X},\mathbf{Y}) = 1$, then there exist a vector $a$, a nonzero real number $b$ and an orthogonal matrix $C$ such that $\mathbf{Y} = a + b\mathbf{X}C$.*

PROOF. In (i), $\mathcal{R}(X,Y)$ exists whenever $X$ and $Y$ have finite first moments, and $X$ and $Y$ are independent if and only if the numerator
$$\mathcal{V}^2(X,Y) = \|f_{X,Y}(t,s) - f_X(t)f_Y(s)\|^2$$
of $\mathcal{R}^2(X,Y)$ is zero. Let $U = e^{i\langle t,X\rangle} - f_X(t)$ and $V = e^{i\langle s,Y\rangle} - f_Y(s)$. Then
$$|f_{X,Y}(t,s) - f_X(t)f_Y(s)|^2 = |E[UV]|^2 \leq (E[|U||V|])^2 \leq E[|U|^2|V|^2]$$
$$= (1 - |f_X(t)|^2)(1 - |f_Y(s)|^2).$$
Thus
$$\int_{\mathbb{R}^{p+q}} |f_{X,Y}(t,s) - f_X(t)f_Y(s)|^2 \, d\omega$$
$$\leq \int_{\mathbb{R}^{p+q}} |(1 - |f_X(t)|^2)(1 - |f_Y(s)|^2)|^2 \, d\omega;$$
hence $0 \leq \mathcal{R}(X,Y) \leq 1$, and (ii) follows by a similar argument.

(iii) If $\mathcal{R}_n(\mathbf{X},\mathbf{Y}) = 1$, then the arguments below show that $X$ and $Y$ are similar almost surely, thus the dimensions of the linear subspaces spanned by $\mathbf{X}$ and $\mathbf{Y}$ respectively are almost surely equal. (Here *similar* means that $\mathbf{Y}$ and $\varepsilon\mathbf{X}$ are isometric for some $\varepsilon \neq 0$.) For simplicity we can suppose that $\mathbf{X}$ and $\mathbf{Y}$ are in the same Euclidean space and both span $\mathbb{R}^p$. From the Cauchy–Bunyakovski inequality it is easy to see that $\mathcal{R}_n(\mathbf{X},\mathbf{Y}) = 1$ if and only if $A_{kl} = \varepsilon B_{kl}$ for some factor $\varepsilon$. Suppose that $|\varepsilon| = 1$. Then
$$|X_k - X_l|_p = |Y_k - Y_l|_q + d_k + d_l$$

for all $k, l$, for some constants $d_k, d_l$. Then with $k = l$ we obtain $d_k = 0$ for all $k$. Now, one can apply a geometric argument. The two samples are isometric, so $\mathbf{Y}$ can be obtained from $\mathbf{X}$ through operations of shift, rotation and reflection, and hence $\mathbf{Y} = a + b\mathbf{X}C$ for some vector $a$, $b = \varepsilon$ and orthogonal matrix $C$. If $|\varepsilon| \neq 1$ and $\varepsilon \neq 0$, apply the geometric argument to $\varepsilon \mathbf{X}$ and $\mathbf{Y}$ and it follows that $\mathbf{Y} = a + b\mathbf{X}C$ where $b = \varepsilon$. □

THEOREM 4 (Properties of dVar). *The following properties hold for random vectors with finite first moments:*

(i) $\mathrm{dVar}(X) = 0$ *implies that* $X = E[X]$, *almost surely.*
(ii) $\mathrm{dVar}(a + bCX) = |b|\,\mathrm{dVar}(X)$ *for all constant vectors $a$ in $\mathbb{R}^p$, scalars $b$ and $p \times p$ orthonormal matrices $C$.*
(iii) $\mathrm{dVar}(X + Y) \leq \mathrm{dVar}(X) + \mathrm{dVar}(Y)$ *for independent random vectors $X$ in $\mathbb{R}^p$ and $Y$ in $\mathbb{R}^p$.*

PROOF. (i) If $\mathrm{dVar}(X) = 0$, then
$$(\mathrm{dVar}(X))^2 = \int_{\mathbb{R}^p} \int_{\mathbb{R}^p} \frac{|f_{X,X}(t,s) - f_X(t)f_X(s)|^2}{c_p^2 |t|_p^{p+1} |s|_p^{p+1}}\, dt\, ds = 0,$$
or equivalently, $f_{X,X}(t,s) = f_X(t+s) = f_X(t)f_X(s)$ for all $t$, $s$. That is, $f_X(t) = e^{i\langle c,t\rangle}$ for some constant vector $c$, and hence $X$ is a constant vector, almost surely.

Statement (ii) is obvious.

(iii) Let $\Delta = f_{X+Y,X+Y}(t,s) - f_X(t+s)f_Y(t+s)$, $\Delta_1 = f_X(t,s) - f_X(t)f_X(s)$ and $\Delta_2 = f_Y(t,s) - f_Y(t)f_Y(s)$. If $X$ and $Y$ are independent, then
$$\begin{aligned}
\Delta &= f_{X+Y,X+Y}(t,s) - f_{X+Y}(t)f_{X+Y}(s) \\
&= f_X(t+s)f_Y(t+s) - f_X(t)f_X(s)f_Y(t)f_Y(s) \\
&= [f_X(t+s) - f_X(t)f_X(s)]f_Y(t+s) \\
&\quad + f_X(t)f_X(s)[f_Y(t+s) - f_Y(t)f_Y(s)] \\
&= \Delta_1 f_Y(t+s) + f_X(t)f_X(s)\Delta_2,
\end{aligned}$$
and therefore $|\Delta|^2 \leq |\Delta_1|^2 + |\Delta_2|^2 + 2|\Delta_1||\Delta_2|$. Equality holds if and only if $\Delta_1 \Delta_2 = 0$, that is, if and only if $X$ or $Y$ is a constant almost surely. □

2.4. *Asymptotic properties of $n\mathcal{V}_n^2$.* Our proposed test of independence is based on the statistic $n\mathcal{V}_n^2/S_2$. If $E(|X|_p + |Y|_q) < \infty$, we prove that under independence $n\mathcal{V}_n^2/S_2$ converges in distribution to a quadratic form

(2.25) $$Q \stackrel{D}{=} \sum_{j=1}^{\infty} \lambda_j Z_j^2,$$

where $Z_j$ are independent standard normal random variables, $\{\lambda_j\}$ are nonnegative constants that depend on the distribution of $(X,Y)$ and $E[Q] = 1$. A test of independence that rejects independence for large $n\mathcal{V}_n^2/S_2$ is statistically consistent against all alternatives with finite first moments.

Let $\zeta(\cdot)$ denote a complex-valued zero-mean Gaussian random process with covariance function

$$R(u, u_0) = (f_X(t - t_0) - f_X(t)\overline{f_X(t_0)})(f_Y(s - s_0) - f_Y(s)\overline{f_Y(s_0)}),$$

where $u = (t, s), u_0 = (t_0, s_0) \in \mathbb{R}^p \times \mathbb{R}^q$.

THEOREM 5 (Weak convergence). *If $X$ and $Y$ are independent and $E(|X|_p + |Y|_q) < \infty$, then*

$$n\mathcal{V}_n^2 \xrightarrow[n\to\infty]{D} \|\zeta(t, s)\|^2.$$

PROOF. Define the empirical process

$$\zeta_n(u) = \zeta_n(t, s) = \sqrt{n}\xi_n(t, s) = \sqrt{n}(f_{X,Y}^n(t, s) - f_X^n(t)f_Y^n(s)).$$

Under the independence hypothesis, $E[\zeta_n(u)] = 0$ and $E[\zeta_n(u)\overline{\zeta_n(u_0)}] = \frac{n-1}{n} \times R(u, u_0)$. In particular, $E|\zeta_n(u)|^2 = \frac{n-1}{n}(1 - |f_X(t)|^2)(1 - |f_Y(s)|^2) \leq 1$.

For each $\delta > 0$ we construct a sequence of random variables $\{Q_n(\delta)\}$ with the following properties:

(i) $Q_n(\delta)$ converges in distribution to a random variable $Q(\delta)$.
(ii) $E|Q_n(\delta) - \zeta_n| \leq \delta$.
(iii) $E|Q(\delta) - \zeta| \leq \delta$.

Then the weak convergence of $\|\zeta_n\|^2$ to $\|\zeta\|^2$ follows from the convergence of the corresponding characteristic functions.

The sequence $Q_n(\delta)$ is defined as follows. Given $\epsilon > 0$, choose a partition $\{D_k\}_{k=1}^N$ of $D(\delta)$ (2.20) into $N = N(\epsilon)$ measurable sets with diameter at most $\epsilon$. Define

$$Q_n(\delta) = \sum_{k=1}^N \int_{D_k} |\zeta_n|^2 \, d\omega.$$

For a fixed $M > 0$ let

$$\beta(\epsilon) = \sup_{u, u_0} E||\zeta_n(u)|^2 - |\zeta_n(u_0)|^2|,$$

where the supremum is taken over all $u = (t, s)$ and $u_0 = (t_0, s_0)$ such that $\max\{|t|, |t_0|, |s|, |s_0|\} < M$, and $|t - t_0|^2 + |s - s_0|^2 < \epsilon^2$. Then $\lim_{\epsilon \to 0} \beta(\epsilon) = 0$ for every fixed $M > 0$, and for fixed $\delta > 0$

$$E\left|\int_{D(\delta)} |\zeta_n(u)|^2 \, d\omega - Q_n(\delta)\right| \leq \beta(\epsilon) \int_D |\zeta_n(u)|^2 \, d\omega \xrightarrow[\epsilon \to 0]{} = 0.$$

On the other hand,

$$\left| \int_D |\zeta_n(u)|^2 \, d\omega - \int_{\mathbb{R}^{p+q}} |\zeta_n(u)|^2 \, d\omega \right|$$
$$\leq \int_{|t|<\delta} |\zeta_n(u)|^2 \, d\omega + \int_{|t|>1/\delta} |\zeta_n(u)|^2 \, d\omega$$
$$+ \int_{|s|<\delta} |\zeta_n(u)|^2 \, d\omega + \int_{|s|>1/\delta} |\zeta_n(u)|^2 \, d\omega.$$

By similar steps as in the proof of Theorem 2, one can derive that

$$E\left[ \int_{|t|<\delta} |\zeta_n(u)|^2 \, d\omega + \int_{|t|>1/\delta} |\zeta_n(u)|^2 \, d\omega \right]$$
$$\leq \frac{n-1}{n}(E|X_1 - X_2| G(|X_1 - X_2|\delta) + w_p \delta) E|Y_1 - Y_2| \xrightarrow[\delta \to 0]{} 0,$$

where $w_p$ is a constant depending only on $p$, and similarly

$$E\left[ \int_{|s|<\delta} |\zeta_n(u)|^2 \, d\omega + \int_{|s|>1/\delta} |\zeta_n(u)|^2 \, d\omega \right] \xrightarrow[\delta \to 0]{} 0.$$

Similar inequalities also hold for the random process $\zeta(t,s)$ with

$$Q(\delta) = \sum_{k=1}^N \int_{D_k} |\zeta(u)|^2 \, d\omega.$$

The weak convergence of $Q_n(\delta)$ to $Q(\delta)$ as $n \to \infty$ follows from the multivariate central limit theorem, and therefore $n\mathcal{V}_n^2 = \|\zeta_n\|^2 \xrightarrow[n \to \infty]{D} \|\zeta\|^2$. □

COROLLARY 2. *If $E(|X|_p + |Y|_q) < \infty$, then:*

(i) *If $X$ and $Y$ are independent, $n\mathcal{V}_n^2/S_2 \xrightarrow[n \to \infty]{D} Q$ where $Q$ is a nonnegative quadratic form of centered Gaussian random variables (2.25) and $E[Q] = 1$.*

(ii) *If $X$ and $Y$ are dependent, then $n\mathcal{V}_n^2/S_2 \xrightarrow[n \to \infty]{P} \infty$.*

PROOF. (i) The independence of $X$ and $Y$ implies that $\zeta_n$ and thus $\zeta$ is a zero-mean process. According to Kuo [5], Chapter 1, Section 2, the squared norm $\|\zeta\|^2$ of the zero-mean Gaussian process $\zeta$ has the representation

$$\|\zeta\|^2 \stackrel{D}{=} \sum_{j=1}^\infty \lambda_j Z_j^2, \tag{2.26}$$

where $Z_j$ are independent standard normal random variables, and the nonnegative constants $\{\lambda_j\}$ depend on the distribution of $(X,Y)$. Hence, under independence, $n\mathcal{V}_n^2$ converges in distribution to a quadratic form (2.26). It follows from (2.5) that

$$E\|\zeta\|^2 = \int_{\mathbb{R}^{p+q}} R(u,u)\,d\omega$$
$$= \int_{\mathbb{R}^{p+q}} (1-|f_X(t)|^2)(1-|f_Y(s)|^2)\,d\omega$$
$$= E(|X-X'|_p|Y-Y'|_q).$$

By the SLLN for $V$-statistics, $S_2 \xrightarrow[n\to\infty]{a.s.} E(|X-X'|_p|Y-Y'|_q)$. Therefore $n\mathcal{V}_n^2/S_2 \xrightarrow[n\to\infty]{D} Q$, where $E[Q] = 1$ and $Q$ is the quadratic form (2.25).

(ii) Suppose that $X$ and $Y$ are dependent and $E(|X|_p + |Y|_q) < \infty$. Then $\mathcal{V}(X,Y) > 0$, Theorem 2 implies that $\mathcal{V}_n^2(\mathbf{X},\mathbf{Y}) \xrightarrow[n\to\infty]{a.s.} \mathcal{V}^2(X,Y) > 0$, and therefore $n\mathcal{V}_n^2(\mathbf{X},\mathbf{Y}) \xrightarrow[n\to\infty]{P} \infty$. By the SLLN, $S_2$ converges to a constant and therefore $n\mathcal{V}_n^2/S_2 \xrightarrow[n\to\infty]{P} \infty$. □

THEOREM 6. *Suppose $T(X,Y,\alpha,n)$ is the test that rejects independence if*

$$\frac{n\mathcal{V}_n^2(\mathbf{X},\mathbf{Y})}{S_2} > (\Phi^{-1}(1-\alpha/2))^2,$$

*where $\Phi(\cdot)$ denotes the standard normal cumulative distribution function, and let $\alpha(X,Y,n)$ denote the achieved significance level of $T(X,Y,\alpha,n)$. If $E(|X|_p + |Y|_q) < \infty$, then for all $0 < \alpha \leq 0.215$*

(i) $\lim_{n\to\infty} \alpha(X,Y,n) \leq \alpha$,
(ii) $\sup_{X,Y}\{\lim_{n\to\infty} \alpha(X,Y,n) : \mathcal{V}(X,Y) = 0\} = \alpha$.

PROOF. (i) The following inequality is proved as a special case of a theorem of Székely and Bakirov [9], page 189. If $Q$ is a quadratic form of centered Gaussian random variables and $E[Q] = 1$, then

$$P\{Q \geq (\Phi^{-1}(1-\alpha/2))^2\} \leq \alpha$$

for all $0 < \alpha \leq 0.215$.

(ii) For Bernoulli random variables $X$ and $Y$ we have that $\mathcal{R}_n(\mathbf{X},\mathbf{Y}) = |\hat{\rho}(\mathbf{X},\mathbf{Y})|$. By the central limit theorem, under independence $\sqrt{n}\hat{\rho}(\mathbf{X},\mathbf{Y})$ is asymptotically normal. Thus, in case $X$ and $Y$ are independent Bernoulli variables, the quadratic form $Q$ contains only one term, $Q = Z_1^2$, and the upper bound $\alpha$ is achieved. □

Thus, a test rejecting independence of $X$ and $Y$ when $\sqrt{n\mathcal{V}_n^2/S_2} \geq \Phi^{-1}(1-\alpha/2)$ has an asymptotic significance level at most $\alpha$. The asymptotic test criterion could be quite conservative for many distributions. Alternatively one can estimate the critical value for the test by conditioning on the observed sample, which is discussed in Section 5.

REMARK 3.  If $E|X|_p^2 < \infty$ and $E|Y|_q^2 < \infty$, then $E[|X|_p|Y|_q] < \infty$, so by Lemma 1 and by Fubini's theorem we can evaluate

$$\mathcal{V}^2(X,Y) = E[|X_1 - X_2|_p |Y_1 - Y_2|_q] + E|X_1 - X_2|_p \, E|Y_1 - Y_2|_q$$
$$- 2E[|X_1 - X_2|_p |Y_1 - Y_3|_q].$$

If second moments exist, Theorem 2 and weak convergence can be established by $V$-statistic limit theorems [13]. Under the null hypothesis of independence, $\mathcal{V}_n^2$ is a degenerate kernel $V$-statistic. The first-order degeneracy follows from inequalities proved in [10]. Thus $n\mathcal{V}_n^2$ converges in distribution to a quadratic form (2.26).

## 3. Extensions.

3.1. *The class of $\alpha$-distance dependence measures.*  We introduce a one-parameter family of distance dependence measures indexed by a positive exponent $\alpha$. In our definition of dCor we have applied exponent $\alpha = 1$.

Suppose that $E(|X|_p^\alpha + |Y|_q^\alpha) < \infty$. Let $\mathcal{V}^{(\alpha)}$ denote the $\alpha$-distance covariance, which is the nonnegative number defined by

$$\mathcal{V}^{2(\alpha)}(X,Y) = \|f_{X,Y}(t,s) - f_X(t)f_Y(s)\|_\alpha^2$$
$$= \frac{1}{C(p,\alpha)C(q,\alpha)} \int_{\mathbb{R}^{p+q}} \frac{|f_{X,Y}(t,s) - f_X(t)f_Y(s)|^2}{|t|_p^{\alpha+p}|s|_q^{\alpha+q}} \, dt\, ds.$$

Similarly, $\mathcal{R}^{(\alpha)}$ denotes $\alpha$-distance correlation, which is the square root of

$$\mathcal{R}^{2(\alpha)} = \frac{\mathcal{V}^{2(\alpha)}(X,Y)}{\sqrt{\mathcal{V}^{2(\alpha)}(X)\mathcal{V}^{2(\alpha)}(Y)}}, \qquad 0 < \mathcal{V}^{2(\alpha)}(X), \qquad \mathcal{V}^{2(\alpha)}(Y) < \infty,$$

and $\mathcal{R}^{(\alpha)} = 0$ if $\mathcal{V}^{2(\alpha)}(X)\mathcal{V}^{2(\alpha)}(Y) = 0$.

The $\alpha$-distance dependence statistics are defined by replacing the exponent 1 with exponent $\alpha$ in the distance dependence statistics (2.8), (2.9) and (2.10). That is, replace $a_{kl} = |X_k - X_l|_p$ with $a_{kl} = |X_k - X_l|_p^\alpha$ and replace $b_{kl} = |Y_k - Y_l|_q$ with $b_{kl} = |Y_k - Y_l|_q^\alpha$, $k, l = 1, \ldots, n$.

Theorem 2 can be generalized for $\|\cdot\|_\alpha$-norms, so that almost sure convergence of $\mathcal{V}_n^{(\alpha)} \to \mathcal{V}^{(\alpha)}$ follows if the $\alpha$-moments are finite. Similarly one can prove the weak convergence and statistical consistency for $\alpha$ exponents, $0 < \alpha < 2$, provided that $\alpha$ moments are finite.

The case $\alpha = 2$ leads to the counterpart of classical correlation and covariance. In fact, if $p = q = 1$, then $\mathcal{R}^{(2)} = |\rho|$, $\mathcal{R}_n^{(2)} = |\hat{\rho}|$ and $\mathcal{V}_n^{(2)} = 2|\hat{\sigma}_{xy}|$, where $\hat{\sigma}_{xy}$ is the maximum likelihood estimator of $\mathrm{Cov}(X,Y)$.

3.2. *Affine invariance.* Group invariance is an important concept in statistical inference (see Eaton [3] or Giri [4]), particularly when any transformation of data and/or parameters by some group element constitutes an equivalent problem for inference. For the problem of testing independence, which is preserved under the group of affine transformations, it is natural to consider dependence measures that are affine invariant. Although $\mathcal{R}(X,Y)$ as defined by (2.7) is not affine invariant, it is clearly invariant with respect to the group of orthogonal transformations

(3.1) $$X \mapsto a_1 + b_1 C_1 X, \qquad Y \mapsto a_2 + b_2 C_2 Y,$$

where $a_1$, $a_2$ are arbitrary vectors, $b_1$, $b_2$ are arbitrary nonzero numbers and $C_1$, $C_2$ are arbitrary orthogonal matrices. We can also define a distance correlation that is affine invariant.

For random samples $\mathbf{X}$ from the distribution of $X$ in $\mathbb{R}^p$ and $\mathbf{Y}$ from the distribution of $Y$ in $\mathbb{R}^q$, define the scaled samples $\mathbf{X}^*$ and $\mathbf{Y}^*$ by

(3.2) $$\mathbf{X}^* = \mathbf{X} S_X^{-1/2}, \qquad \mathbf{Y}^* = \mathbf{Y} S_Y^{-1/2},$$

where $S_X$ and $S_Y$ are the sample covariance matrices of $\mathbf{X}$ and $\mathbf{Y}$, respectively. Although the sample vectors in (3.2) are not invariant to affine transformations, the distances $|X_k^* - X_l^*|$ and $|Y_k^* - Y_l^*|$, $k, l = 1, \ldots, n$, are invariant to affine transformations. Then the affine distance correlation statistic $\mathcal{R}_n^*(\mathbf{X}, \mathbf{Y})$ between random samples $\mathbf{X}$ and $\mathbf{Y}$ is the square root of

$$\mathcal{R}_n^{*2}(\mathbf{X}, \mathbf{Y}) = \frac{\mathcal{V}_n^2(\mathbf{X}^*, \mathbf{Y}^*)}{\sqrt{\mathcal{V}_n^2(\mathbf{X}^*) \mathcal{V}_n^2(\mathbf{Y}^*)}}.$$

Properties established in Section 2 also hold for $\mathcal{V}_n^*$ and $\mathcal{R}_n^*$, because the transformation (3.2) simply replaces the weight function $\{c_p c_q |t|_p^{1+p} |s|_q^{1+q}\}^{-1}$ with the weight function $\{c_p c_q |S_X^{1/2} t|_p^{1+p} |S_Y^{1/2} s|_q^{1+q}\}^{-1}$.

**4. Results for the bivariate normal distribution.** Let $X$ and $Y$ have standard normal distributions with $\mathrm{Cov}(X,Y) = \rho(X,Y) = \rho$. Introduce the function

$$F(\rho) = \int_{-\infty}^{\infty} \int_{-\infty}^{\infty} |f_{X,Y}(t,s) - f_X(t) f_Y(s)|^2 \frac{dt}{t^2} \frac{ds}{s^2}.$$

Then $\mathcal{V}^2(X,Y) = F(\rho)/c_1^2 = F(\rho)/\pi^2$ and

(4.1) $$\mathcal{R}^2(X,Y) = \frac{\mathcal{V}^2(X,Y)}{\sqrt{\mathcal{V}^2(X,X) \mathcal{V}^2(Y,Y)}} = \frac{F(\rho)}{F(1)}.$$

THEOREM 7. *If $X$ and $Y$ are standard normal, with correlation $\rho = \rho(X,Y)$, then*

(i) $\mathcal{R}(X,Y) \leq |\rho|$,
(ii) $\mathcal{R}^2(X,Y) = \frac{\rho \arcsin \rho + \sqrt{1-\rho^2} - \rho \arcsin \rho/2 - \sqrt{4-\rho^2} + 1}{1 + \pi/3 - \sqrt{3}}$,
(iii) $\inf_{\rho \neq 0} \frac{\mathcal{R}(X,Y)}{|\rho|} = \lim_{\rho \to 0} \frac{\mathcal{R}(X,Y)}{|\rho|} = \frac{1}{2(1+\pi/3-\sqrt{3})^{1/2}} \approx 0.89066$.

PROOF. (i) If $X$ and $Y$ are standard normal with correlation $\rho$, then

$$F(\rho) = \int_{-\infty}^{\infty} \int_{-\infty}^{\infty} |e^{-(t^2+s^2)/2-\rho ts} - e^{-t^2/2} e^{-s^2/2}|^2 \frac{dt}{t^2} \frac{ds}{s^2}$$

$$= \int_{\mathbb{R}^2} e^{-t^2-s^2}(1 - 2e^{-\rho ts} + e^{-2\rho ts}) \frac{dt}{t^2} \frac{ds}{s^2}$$

$$= \int_{\mathbb{R}^2} e^{-t^2-s^2} \sum_{n=2}^{\infty} \frac{2^n - 2}{n!} (-\rho ts)^n \frac{dt}{t^2} \frac{ds}{s^2}$$

$$= \int_{\mathbb{R}^2} e^{-t^2-s^2} \sum_{k=1}^{\infty} \frac{2^{2k} - 2}{(2k)!} (-\rho ts)^{2k} \frac{dt}{t^2} \frac{ds}{s^2}$$

$$= \rho^2 \left[ \sum_{k=1}^{\infty} \frac{2^{2k} - 2}{(2k)!} \rho^{2(k-1)} \int_{\mathbb{R}^2} e^{-t^2-s^2} (ts)^{2(k-1)} \, dt \, ds \right].$$

Thus $F(\rho) = \rho^2 G(\rho)$, where $G(\rho)$ is a sum with all nonnegative terms. The function $G(\rho)$ is clearly nondecreasing in $\rho$ and $G(\rho) \leq G(1)$. Therefore

$$\mathcal{R}^2(X,Y) = \frac{F(\rho)}{F(1)} = \rho^2 \frac{G(\rho)}{G(1)} \leq \rho^2,$$

or equivalently, $\mathcal{R}(X,Y) \leq |\rho|$.

(ii) Note that $F(0) = F'(0) = 0$ so $F(\rho) = \int_0^\rho \int_0^x F''(z) \, dz \, dx$. The second derivative of $F$ is

$$F''(z) = \frac{d^2}{dz^2} \int_{\mathbb{R}^2} e^{-t^2-s^2}(1 - 2e^{-zts} + e^{-2zts}) \frac{dt}{t^2} \frac{ds}{s^2} = 4V(z) - 2V\left(\frac{z}{2}\right),$$

where

$$V(z) = \int_{\mathbb{R}^2} e^{-t^2-s^2-2zts} \, dt \, ds = \frac{\pi}{\sqrt{1-z^2}}.$$

Here we have applied a change of variables, used the fact that the eigenvalues of the quadratic form $t^2 + s^2 + 2zts$ are $1 \pm z$, and $\int_{-\infty}^{\infty} e^{-t^2 \lambda} \, dt = (\pi/\lambda)^{1/2}$. Then

$$F(\rho) = \int_0^\rho \int_0^x \left( \frac{4\pi}{\sqrt{1-z^2}} - \frac{2\pi}{\sqrt{1-z^2/4}} \right) dz \, dx$$

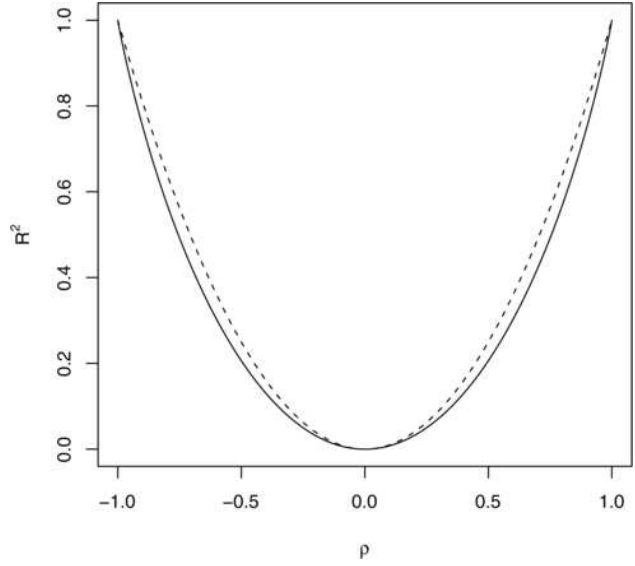

Fig. 1. *Dependence coefficient $\mathcal{R}^2$ (solid line) and correlation $\rho^2$ (dashed line) in the bivariate normal case.*

$$= 4\pi \int_0^\rho (\arcsin(x) - \arcsin(x/2))\,dx$$

$$= 4\pi(\rho \arcsin \rho + \sqrt{1-\rho^2} - \rho \arcsin(\rho/2) - \sqrt{4-\rho^2} + 1),$$

and (4.1) implies (ii).

(iii) In the proof of (i) we have that $\mathcal{R}/|\rho|$ is a nondecreasing function of $|\rho|$, and $\lim_{|\rho| \to 0} \mathcal{R}(X,Y)/|\rho| = (1 + \pi/3 - \sqrt{3})^{-1/2}/2$ follows from (ii). □

The relation between $\mathcal{R}$ and $\rho$ derived in Theorem 7 is shown by the plot of $\mathcal{R}^2$ versus $\rho^2$ in Figure 1.

**5. Empirical results.** In this section we summarize Monte Carlo power comparisons of our proposed distance covariance test with three classical tests for multivariate independence. The likelihood ratio test (LRT) of the hypothesis $H_0: \Sigma_{12} = 0$, with $\mu$ unknown, is based on

$$(5.1) \qquad \frac{\det(S)}{\det(S_{11})\det(S_{22})} = \frac{\det(S_{22} - S_{21}S_{11}^{-1}S_{12})}{\det(S_{22})},$$

where $\det(\cdot)$ is the determinant, $S$, $S_{11}$ and $S_{22}$ denote the sample covariances of $(\mathbf{X},\mathbf{Y})$, $\mathbf{X}$ and $\mathbf{Y}$, respectively, and $S_{12}$ is the sample covariance $\widehat{\mathrm{Cov}}(\mathbf{X},\mathbf{Y})$. Under multivariate normality,

$$W = 2\log \lambda = -n \log \det(I - S_{22}^{-1}S_{21}S_{11}^{-1}S_{12})$$

TABLE 1
*Empirical Type-*I *error rates for 10,000 tests at nominal significance level* 0.1 *in Example 1 ($\rho = 0$), using B replicates for V and Bartlett's approximation for W, T and S*

|       |     | (a) Multivariate normal, $p = q = 5$ | | | | (b) $t(1)$, $p = q = 5$ | | | |
|-------|-----|--------|--------|--------|--------|--------|--------|--------|--------|
| $n$   | $B$ | $V$    | $W$    | $T$    | $S$    | $V$    | $W$    | $T$    | $S$    |
| 25    | 400 | 0.1039 | 0.1089 | 0.1212 | 0.1121 | 0.1010 | 0.3148 | 0.1137 | 0.1120 |
| 30    | 366 | 0.0992 | 0.0987 | 0.1145 | 0.1049 | 0.0984 | 0.3097 | 0.1102 | 0.1078 |
| 35    | 342 | 0.0977 | 0.1038 | 0.1091 | 0.1011 | 0.1060 | 0.3102 | 0.1087 | 0.1054 |
| 50    | 300 | 0.0990 | 0.0953 | 0.1052 | 0.1011 | 0.1036 | 0.2904 | 0.1072 | 0.1037 |
| 70    | 271 | 0.1001 | 0.0983 | 0.1031 | 0.1004 | 0.1000 | 0.2662 | 0.1013 | 0.0980 |
| 100   | 250 | 0.1019 | 0.0954 | 0.0985 | 0.0972 | 0.1025 | 0.2433 | 0.0974 | 0.1019 |
|       |     | (c) $t(2)$, $p = q = 5$ | | | | (d) $t(3)$, $p = q = 5$ | | | |
| $n$   | $B$ | $V$    | $W$    | $T$    | $S$    | $V$    | $W$    | $T$    | $S$    |
| 25    | 400 | 0.1037 | 0.1612 | 0.1217 | 0.1203 | 0.1001 | 0.1220 | 0.1174 | 0.1134 |
| 30    | 366 | 0.0959 | 0.1636 | 0.1070 | 0.1060 | 0.1017 | 0.1224 | 0.1143 | 0.1062 |
| 35    | 342 | 0.0998 | 0.1618 | 0.1080 | 0.1073 | 0.1074 | 0.1213 | 0.1131 | 0.1047 |
| 50    | 300 | 0.1033 | 0.1639 | 0.1010 | 0.0969 | 0.1050 | 0.1166 | 0.1065 | 0.1046 |
| 70    | 271 | 0.1029 | 0.1590 | 0.1063 | 0.0994 | 0.1037 | 0.1200 | 0.1020 | 0.1002 |
| 100   | 250 | 0.0985 | 0.1560 | 0.1050 | 0.1007 | 0.1019 | 0.1176 | 0.1066 | 0.1033 |

has the Wilks Lambda distribution $\Lambda(q, n - 1 - p, p)$ [14]. Puri and Sen [8], Chapter 8, proposed similar tests based on more general sample dispersion matrices $T = (T_{ij})$. The Puri–Sen tests replace $S$, $S_{11}$, $S_{12}$ and $S_{22}$ in (5.1) with $T$, $T_{11}$, $T_{12}$ and $T_{22}$. For example, $T$ can be a matrix of Spearman's rank correlation statistics. For a sign test the dispersion matrix has entries $\frac{1}{n}\sum_{j=1}^{n}\mathrm{sign}(Z_{jk} - \widetilde{Z}_k)\mathrm{sign}(Z_{jm} - \widetilde{Z}_m)$, where $\widetilde{Z}_k$ is the sample median of the $k$th variable. Critical values of the Wilks Lambda and Puri–Sen statistics are given by Bartlett's approximation: if $n$ is large and $p, q > 2$, then $-(n - \frac{1}{2}(p + q + 3))\log\det(I - S_{22}^{-1}S_{21}S_{11}^{-1}S_{12})$ has an approximate $\chi^2(pq)$ distribution ([6], Section 5.3.2b).

To implement the distance covariance test for small samples, we obtain a reference distribution for $n\mathcal{V}_n^2$ under independence by conditioning on the observed sample, that is, by computing replicates of $n\mathcal{V}_n^2$ under random permutations of the indices of the $Y$ sample. We obtain good control of Type-I error (see Table 1) with a small number of replicates; for this study we used $\lfloor 200 + 5000/n \rfloor$ replicates. Implementation is straightforward due to the simple form of the statistic. The statistic $n\mathcal{V}_n^2$ has $O(n^2)$ time and space computational complexity. Source code for the test implementation is available from the authors upon request.

Each example compares the empirical power of the dCov test (labeled $V$) with the Wilks Lambda statistic ($W$), Puri–Sen rank correlation statistic ($S$) and Puri–Sen sign statistic ($T$). Empirical power is computed as the proportion of significant tests on 10,000 random samples at significance level 0.1.

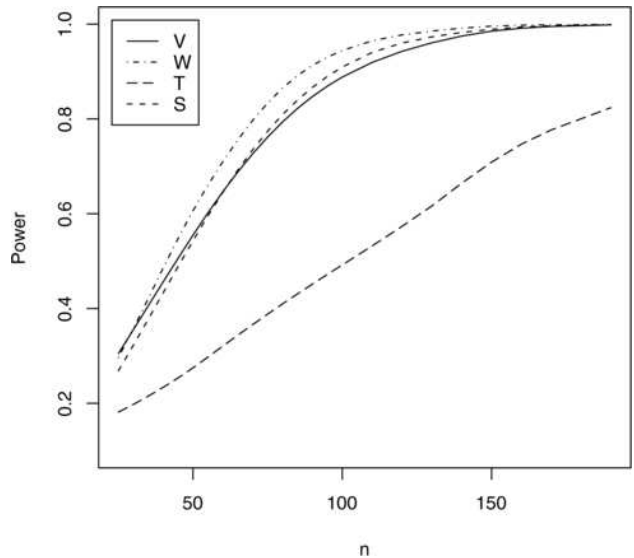

FIG. 2. *Example 1*(a): *Empirical power at* 0.1 *significance and sample size n, multivariate normal alternative.*

EXAMPLE 1. In 1(a) the marginal distributions of $X$ and $Y$ are standard multivariate normal in dimensions $p = q = 5$ and $\text{Cov}(X_k, Y_l) = \rho$ for $k, l = 1, \ldots, 5$. The results displayed in Figure 2 are based on 10,000 tests for each of the sample sizes $n = 25:50:1$, $55:100:5$, $110:200:10$ with $\rho = 0.1$. As expected, the Wilks LRT is optimal in this case, but power of the dCov test is quite close to $W$. Table 1 gives empirical Type-I error rates for this example when $\rho = 0$.

In Examples 1(b)–1(d) we repeat 1(a) under identical conditions except that the random variables $X_k$ and $Y_l$ are generated from the $t(\nu)$ distribution. Table 1 gives empirical Type-I error rates for $\nu = 1, 2, 3$ when $\rho = 0$. Empirical power for the alternative $\rho = 0.1$ is compared in Figures 3–5. (The Wilks LRT has inflated Type-I error for $\nu = 1, 2, 3$, so a power comparison with $W$ is not meaningful, particularly for $\nu = 1, 2$.)

EXAMPLE 2. The distribution of $X$ is standard multivariate normal ($p = 5$), and $Y_{kj} = X_{kj}\varepsilon_{kj}$, $j = 1, \ldots, p$, where $\varepsilon_{kj}$ are independent standard normal variables and independent of $X$. Comparisons are based on 10,000 tests for each of the sample sizes $n = 25:50:1$, $55:100:5$, $110:240:10$. The results displayed in Figure 6 show that the dCov test is clearly superior to the LRT tests. This alternative is an example where the rank correlation and sign tests do not exhibit power increasing with sample size.

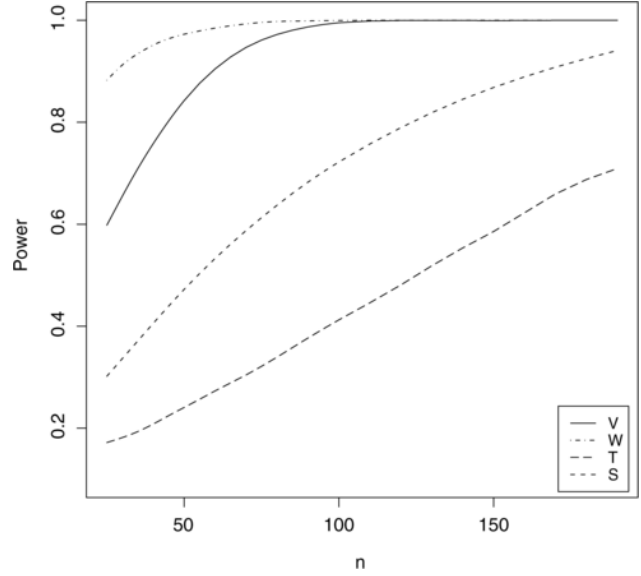

FIG. 3. *Example 1(b): Empirical power at* 0.1 *significance and sample size* $n$, $t(1)$ *alternative.*

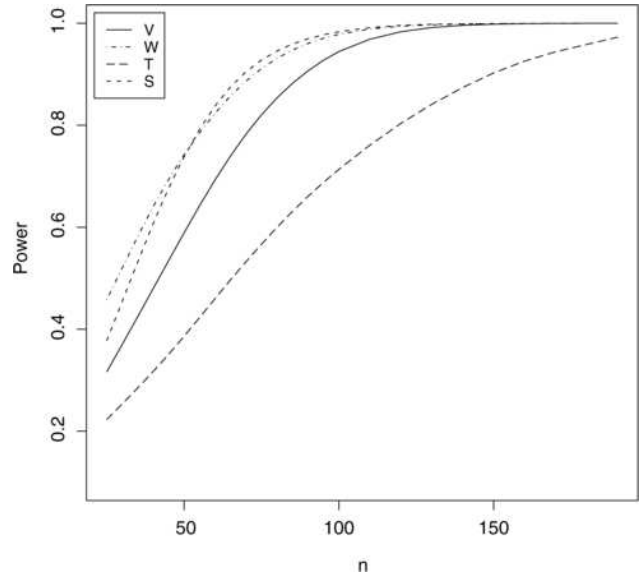

FIG. 4. *Example 1(c): Empirical power at* 0.1 *significance and sample size* $n$, $t(2)$ *alternative.*

EXAMPLE 3. The distribution of $X$ is standard multivariate normal ($p = 5$), and $Y_{kj} = \log(X_{kj}^2)$, $j = 1, \ldots, p$. Comparisons are based on 10,000 tests for each of the sample sizes $n = 25:50:1$, $55:100:5$. Simulation results are displayed in Figure 7. This is an example of a nonlinear relation where $n\mathcal{V}_n^2$ achieves very good power while none of the LRT type tests performs well.

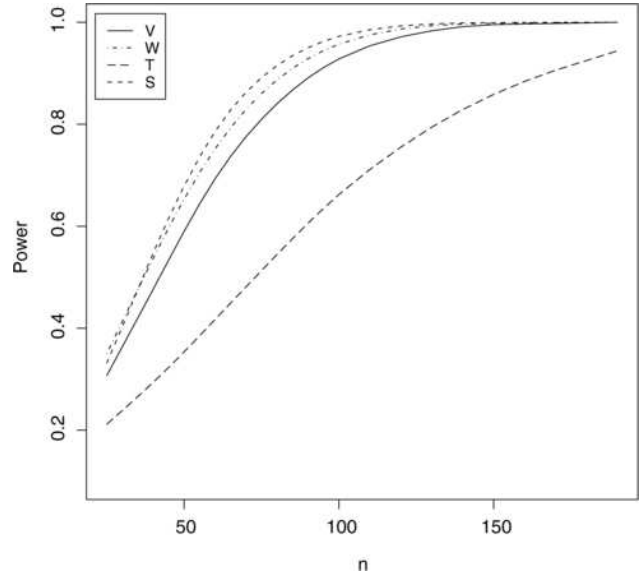

Fig. 5. *Example 1(d): Empirical power at* 0.1 *significance and sample size* $n$, $t(3)$ *alternative.*

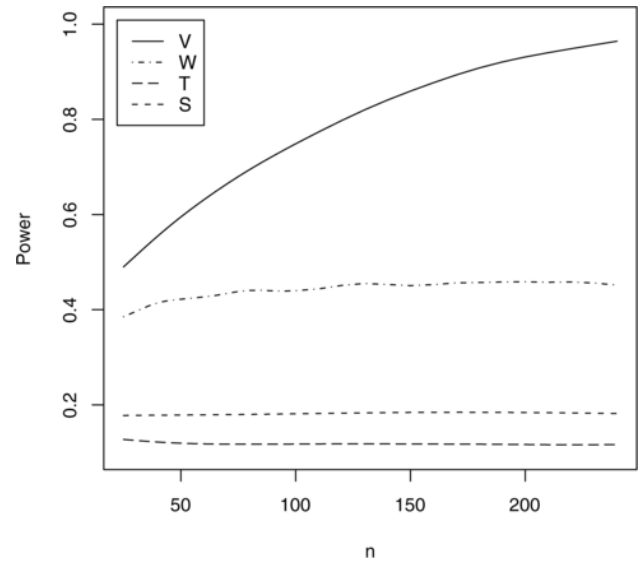

Fig. 6. *Example 2: Empirical power at* 0.1 *significance and sample size* $n$ *against the alternative* $Y = X\varepsilon$.

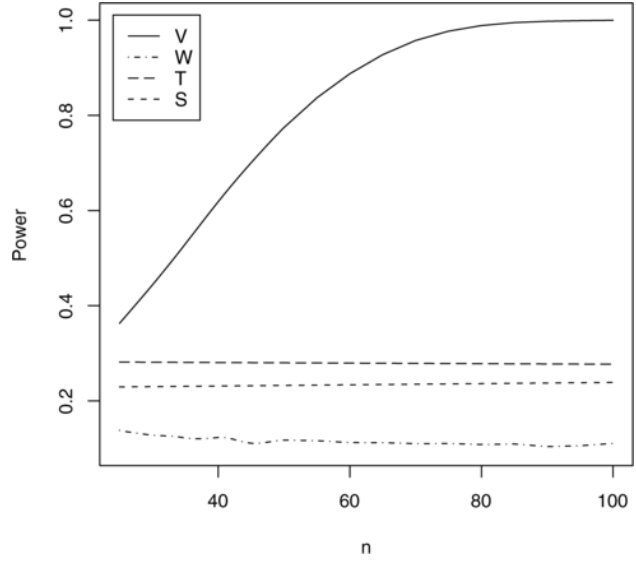

FIG. 7. *Example 3: Empirical power at* 0.1 *significance and sample size n against the alternative* $Y = \log(X^2)$.

**6. Summary.** We have introduced new distance measures of dependence dCov, analogous to covariance, and dCor, analogous to correlation, defined for all random vectors with finite first moments. The dCov test of multivariate independence based on the statistic $n\mathcal{V}_n^2$ is statistically consistent against all dependent alternatives with finite expectation. Empirical results suggest that the dCov test may be more powerful than the parametric LRT when the dependence structure is nonlinear, while in the multivariate normal case the dCov test was quite close in power to the likelihood ratio test. Our proposed statistics are sensitive to all types of departures from independence, including nonlinear or nonmonotone dependence structure.

## APPENDIX

PROOF OF STATEMENT (2.18). By definition (2.8)

$$n^2 \mathcal{V}_n^2 = \sum_{k,l=1}^{n} \left\{ \begin{array}{cccc} a_{kl}b_{kl} & -a_{kl}\bar{b}_{k.} & -a_{kl}\bar{b}_{.l} & +a_{kl}\bar{b}_{..} \\ -\bar{a}_{k.}b_{kl} & +\bar{a}_{k.}\bar{b}_{k.} & +\bar{a}_{k.}\bar{b}_{.l} & -\bar{a}_{k.}\bar{b}_{..} \\ -\bar{a}_{.l}b_{kl} & +\bar{a}_{.l}\bar{b}_{k.} & +\bar{a}_{.l}\bar{b}_{.l} & -\bar{a}_{.l}\bar{b}_{..} \\ +\bar{a}_{..}b_{kl} & -\bar{a}_{..}\bar{b}_{k.} & -\bar{a}_{..}\bar{b}_{.l} & +\bar{a}_{..}\bar{b}_{..} \end{array} \right\}$$

$$= \sum_{k,l} a_{kl}b_{kl} - \sum_{k} a_{k.}\bar{b}_{k.} - \sum_{l} a_{.l}\bar{b}_{.l} + a_{..}\bar{b}_{..}$$

(A.1)
$$- \sum_{k} \bar{a}_{k.}b_{k.} + n \sum_{k} \bar{a}_{k.}\bar{b}_{k.} + \sum_{k,l} \bar{a}_{k.}\bar{b}_{.l} - n \sum_{k} \bar{a}_{k.}\bar{b}_{..}$$

$$- \sum_l \bar{a}_{\cdot l} b_{\cdot l} + \sum_{k,l} \bar{a}_{\cdot l} \bar{b}_{k\cdot} + n \sum_l \bar{a}_{\cdot l} \bar{b}_{\cdot l} - n \sum_l \bar{a}_{\cdot l} \bar{b}_{\cdot \cdot}$$

$$+ \bar{a}_{\cdot\cdot} b_{\cdot\cdot} - n \sum_k \bar{a}_{\cdot\cdot} \bar{b}_{k\cdot} - n \sum_l \bar{a}_{\cdot\cdot} \bar{b}_{\cdot l} + n^2 \bar{a}_{\cdot\cdot} \bar{b}_{\cdot\cdot},$$

where $a_{k\cdot} = n\bar{a}_{k\cdot}$, $a_{\cdot l} = n\bar{a}_{\cdot l}$, $b_{k\cdot} = n\bar{b}_{k\cdot}$ and $b_{\cdot l} = n\bar{b}_{\cdot l}$. Applying the identities

$$(A.2) \qquad S_1 = \frac{1}{n^2} \sum_{k,l=1}^n a_{kl} b_{kl},$$

$$(A.3) \qquad S_2 = \frac{1}{n^2} \sum_{k,l=1}^n a_{kl} \frac{1}{n^2} \sum_{k,l=1}^n b_{kl} = \bar{a}_{\cdot\cdot} \bar{b}_{\cdot\cdot},$$

$$(A.4) \qquad n^2 S_2 = n^2 \bar{a}_{\cdot\cdot} \bar{b}_{\cdot\cdot} = \frac{a_{\cdot\cdot}}{n} \sum_{l=1}^n \frac{b_{\cdot l}}{n} = \sum_{k,l=1}^n \bar{a}_{k\cdot} \bar{b}_{\cdot l},$$

$$(A.5) \qquad S_3 = \frac{1}{n^3} \sum_{k=1}^n \sum_{l,m=1}^n a_{kl} b_{km} = \frac{1}{n^3} \sum_{k=1}^n a_{k\cdot} b_{k\cdot} = \frac{1}{n} \sum_{k=1}^n \bar{a}_{k\cdot} \bar{b}_{k\cdot}.$$

to (A.1), we obtain

$$\begin{aligned} n^2 \mathcal{V}_n^2 &= n^2 S_1 - n^2 S_3 - n^2 S_3 + n^2 S_2 \\ &\quad - n^2 S_3 + n^2 S_3 + n^2 S_2 - n^2 S_2 \\ &\quad - n^2 S_3 + n^2 S_2 + n^2 S_3 - n^2 S_2 \\ &\quad + n^2 S_2 - n^2 S_2 - n^2 S_2 + n^2 S_2 = n^2(S_1 + S_2 - 2S_3). \qquad \square \end{aligned}$$

**Acknowledgments.** The authors are grateful to the Associate Editor and referee for many suggestions that greatly improved the paper.

G. J. SZÉKELY  
DEPARTMENT OF MATHEMATICS AND STATISTICS  
BOWLING GREEN STATE UNIVERSITY  
BOWLING GREEN, OHIO 43403  
USA  
AND  
RÉNYI INSTITUTE OF MATHEMATICS  
HUNGARIAN ACADEMY OF SCIENCES  
HUNGARY  
E-MAIL: gabors@bgnet.bgsu.edu

M. L. RIZZO  
DEPARTMENT OF MATHEMATICS AND STATISTICS  
BOWLING GREEN STATE UNIVERSITY  
BOWLING GREEN, OHIO 43403  
USA  
E-MAIL: mrizzo@bgnet.bgsu.edu

N. K. BAKIROV  
INSTITUTE OF MATHEMATICS  
USC RUSSIAN ACADEMY OF SCIENCES  
112 CHERNYSHEVSKII ST.  
450000 UFA  
RUSSIA